\documentclass[journal]{IEEEtran}
\IEEEoverridecommandlockouts 
\usepackage{textcomp}
\usepackage{graphicx} 
\usepackage{subfigure}
\usepackage{wrapfig} 
\usepackage{indentfirst}
\usepackage{amsmath, amssymb,amsfonts}
\usepackage{textcomp}
\usepackage{color}
\usepackage{multicol}
\usepackage{cite}
\usepackage{array}
\usepackage{mathtools}
\usepackage{fancyhdr}
\usepackage{multicol}
\usepackage{nameref}
\usepackage[ruled,vlined]{algorithm2e}
\DeclarePairedDelimiterX{\inp}[2]{\langle}{\rangle}{#1, #2}

\allowdisplaybreaks[3]

\newcommand{\utwi}[1]{\mbox{\boldmath $#1$}}

\renewcommand{\hat}{\widehat}
\renewcommand{\tilde}{\widetilde}

\newcommand{\cT}{{\cal T}}

\newcommand{\cX}{{\cal X}}

\newcommand{\bd}{{\bf d}}

\newcommand{\bs}{{\bf s}}
\newcommand{\bx}{{\bf x}}

\newcommand{\bv}{{\bf v}}
\newcommand{\bw}{{\bf w}}

\newcommand{\bz}{{\bf z}}
\newcommand{\by}{{\bf y}}
\newcommand{\bA}{{\bf A}}
\newcommand{\bB}{{\bf B}}

\newcommand{\bmu}{{\utwi{\mu}}}

\DeclarePairedDelimiterX{\norm}[1]{\lVert}{\rVert}{#1}

\newtheorem{proposition}{Proposition}

\newtheorem{corollary}{Corollary}

\def\BibTeX{{\rm B\kern-.05em{\sc i\kern-.025em b}\kern-.08em
    T\kern-.1667em\lower.7ex\hbox{E}\kern-.125emX}}

\begin{document}

\title{Personalized Demand Response \\via Shape-Constrained Online Learning}
\author{Ana M. Ospina, Andrea Simonetto, and Emiliano Dall'Anese
\thanks{A. Ospina and E. Dall'Anese are with the Department of Electrical, Computer and Energy Engineering, University of Colorado, Boulder, CO, USA; e-mails: ana.ospina, emiliano.dallanese@colorado.edu. A. Simonetto is with IBM Research Ireland, Dublin, Ireland; e-mail: andrea.simonetto@ibm.com.}
\thanks{The work of A. Ospina and E. Dall'Anese was supported by the National Science Foundation CAREER award 1941896 and by the U.S Department of Energy project 3.2.6.80.}
\vspace*{-.5cm} }
\maketitle

\begin{abstract}
This paper formalizes a demand response task as an optimization problem featuring a known time-varying engineering cost and an unknown (dis)comfort function. Based on this model, this paper develops a feedback-based  projected gradient method to solve the 
demand response problem in an online fashion, where: i) feedback from the user is leveraged to learn the (dis)comfort function concurrently with the execution of the algorithm; and, ii) measurements of electrical quantities are used to estimate the gradient of the known engineering cost. To learn the unknown function, a shape-constrained Gaussian Process is leveraged; this approach allows one to obtain an estimated function that is strongly convex and smooth. The performance of the online algorithm is analyzed by using metrics such as the tracking error and the dynamic regret. A numerical example is illustrated to corroborate the technical findings. 


\end{abstract}

%

\section{Introduction}

Net-load and demand response (DR) strategies hold promise to increase the  flexibility and efficiency  of power  systems by allowing controllable devices to provided services at  various time-scales -- from real-time frequency and voltage support to a slower time-scale peak-shifting service \cite{DR2, DR5,Lesage_2020,Bahrami18}. Typical DR formulations involve a composite cost function to strike a balance between system-level operational objectives and (dis)satisfaction of the device's owner~\cite{Gatsis_DR,Poor_2018}; e.g., deviations from a preferred indoor temperature or charging profile of the electric vehicle. This aspect renders the actual implementation of DR programs challenging: users’ preferences, satisfaction and responsiveness to pricing\cite{xu2019learning} are not easy to model; synthetic cost functions adopted in existing demand response and net-load management frameworks favor computational tractability, but may not capture the users’ goals truthfully.
 
In this context, this paper formalizes a DR task as an optimization problem featuring a \emph{known} time-varying engineering cost and an \emph{unknown} (dis)comfort function. The engineering cost can be related to operational efficiency and may capture objectives such as aggregate setpoint tracking when devices aggregate in a virtual power-plant fashion; it is time-varying~\cite{Mag_signal_2020} in a sense that it captures time-varying objectives (e.g., tracking of a power setpoint that evolves over time), dynamic pricing, or real time measurements. In lieu of synthetic mathematical models for the user's functions (based on e.g., statistics or averaged models), this paper leverages Gaussian Processes (GPs) \cite{GPML_Rasmussen06,GP_userfeedback} to learn the function from data (e.g., users' feedback). 
Approximating a function with a GP often leads to a nonconvex smooth cost; to favor computational tractability, and since user's preferences are often well approximated by convex functions (see, e.g.,\cite{Johnson2018} and references therein), we leverage a shape-constrained GP approach where the discomfort function is approximated with a function that is strongly convex, differentiable, and with a Lipschitz gradient~\cite{jrnl_wang}. The paper then develops a feedback-based  projected gradient method to solve the 
demand response problem in an online fashion. The proposed strategy allows to overcome the following challenges:

\noindent \textit{C1)} \textit{Discomfort function uncertainty:} The functions that model the users' discomfort may not be known and models may be inaccurate. Feedback from the user is leveraged to learn the (dis)comfort function concurrently with the execution of the algorithm using a shape-constrained GP.  

\noindent \textit{C2)} \textit{Pervasive metering:} To solve the optimization problem, one may require the measurements of the powers of non-controllable loads at all locations in real time, and this is a problematic task in power systems. In the proposed strategy, measurements of electrical quantities are used to estimate the gradient of the known engineering cost, and information about the non-controllable loads is not necessary. 
 
Examples of related works on real-time DR include the online convex optimization strategy applied to DR problems in \cite{Lesage_2020}; however, the function associated with heating, ventilation, and air-conditioning (HVAC) systems of commercial buildings is known, and no measurements are utilized in the algorithm.  An online learning approach for computing users' optimal scheduling policy were investigated in~\cite{Bahrami18}, for a given householder's cost function. Also, an online learning approach was considered in~\cite{wang2014adaptive}, based on  a
multi-armed restless bandit problem with controlled bandits. Price responsiveness of the end users that participate in DR programs was studied in, e.g., \cite{xu2019learning} by using a dynamical model that captures the temporal behavior of the users. Community-level energy management systems that weakly control consumers were investigated in~\cite{Shibasaki_2019}. 
For completeness, we point out that  users' perception was incorporated in the decision making process with GPs in other application domains as discussed in, e.g.,~\cite{learn1, GP_learn}.

\section{Preliminaries and Problem Statement}

\subsection{Modeling} 
\label{sec:modeling}

We consider a power network with $M$ controllable loads or DERs   -- hereafter refereed to as ``devices'' for brevity\footnote{\textit{Notation:} Upper-case (lower-case) boldface letters will be used for matrices (column vectors), and $(\cdot)^\top$ denotes transposition. For a given column vector $\mathbf{x} \in \mathbb{R}^n$, $\norm{\mathbf{x}} := \sqrt{\mathbf{x}^\top\mathbf{x}}$. A vector of zeros is represented by $\mathbf{0}$ and a vector of ones by $\mathbf{1}$, with the corresponding dimensions. $\mathcal{O}$ refers to the big O notation; that is, given two positive sequences $\{a_k\}_{k = 0}^\infty$ and $\{b_k\}_{k = 0}^\infty$, we say that $a_k = \mathcal{O}(b_k)$ is $\limsup_{k \rightarrow \infty}  (a_k/b_k) < \infty$.}. Time  is discretized as $t \in \cT := \{k \Delta, k \in \mathbb{N} \}$, where $\Delta$ is a given time interval (e.g., one second or a few seconds~\cite{Dallanese2018feedback, Lesage_2020}). Commands are dispatched to the DERs at each time $t$, and the commanded setpoint for the $m$th device is denoted as $x_{m,t} \in \mathcal{X}_{m,t}$, where $\mathcal{X}_{m,t} \subseteq \mathbb{R}$ is a convex and compact set modeling hardware or operational constraints (e.g., real power commands or temperature setpoints). If a device (e.g., a load) can be controlled at the slower rate (e.g., at the minute-level), the respective setpoint is obviously kept constant over a number of time steps (i.e., $\mathcal{X}_{m,t}$ is a singleton set).  To simplify the notation, the setpoints at time $t$ are aggregated in the column vector $\bx_t = [x_{1,t}, x_{2,t}, \dots, x_{M,t}]^\top \in \mathcal{X}_t \subseteq \mathbb{R}^{M}$, where $\mathcal{X}_t$ is a convex and compact set and is defined as $\mathcal{X}_t:= \mathcal{X}_{1,t} \times \mathcal{X}_{2,t} \times \dots \times \mathcal{X}_{M,t}$. 

The setpoints $\mathbf{x}_t$ are mapped to pertinent electrical states $\mathbf{y}_t \in \mathbb{R}^S$ through a mapping $\mathbf{y}_t = \mathcal{M}(\mathbf{x}_t, \mathbf{w}_t)$, where $\mathcal{M}: \mathbb{R}^M \times \mathbb{R}^W \mapsto \mathbb{R}^S$ models the power network effects and $\mathbf{w}_t \in \mathbb{R}^W$ is a (possibly high-dimensional) vector of powers consumed by $W$ non-controllable devices.  In particular, in this paper we focus on a model of the form:
\begin{equation}
\label{eq:constraint}
    \mathbf{y}_t = \mathbf{A} \mathbf{x}_t + \mathbf{B}\mathbf{w}_t,
\end{equation}
where $\mathbf{A} \in \mathbb{R}^{S \times N}$ and $\mathbf{B} \in \mathbb{R}^{S \times W}$ are known (and possibly time-varying) network matrices. Examples for how to build these matrices will be provided shortly. 

The objective is to formulate a demand-side management problem \cite{Gatsis_DR, Poor_2018} that allows real-time scheduling of end-user devices by minimizing a cost that accounts for both network performance metrics and user satisfaction. Accordingly, let $U_m: \mathcal{X}_{m} \mapsto \mathbb{R}$ be a ``discomfort function''  for the the $m$th user or device. For example, for a thermostatically controllable load, this function may model the discomfort of the user for deviations from a preferred setpoint; for an electric vehicle, $U_m$ may model the dissatisfaction of the user for deviations relative to a preferred charging profile\footnote{The function $U_m$ is assumed to be time-invariant for simplicity; however, the proposed approach can be naturally extended to cases where some of the functions $U_{m,t}$ are time-varying functions to model a dynamic user behavior.}.  Many exiting works presume that the function $U_m$ is \emph{known} and it is convex;  as explained shortly, here $U_m$ will be \emph{learned} from data.

Consider the following \emph{time-varying}  problem~\cite{Mag_signal_2020}: 
\begin{subequations}
\label{eq:problem_DR}
\begin{align}
\underset{\{\mathbf{y}_t \in \mathbb{R}^S,\mathbf{x}_t \in \mathcal{X}_t\}_{i=1}^T}{\min} & \quad \sum_{m=1}^M U_m(x_{m,t}) + C_t(\mathbf{y}_t) \\
  \text{subject to: } ~& \mathbf{y}_t =  \mathbf{A} \mathbf{x}_t + \mathbf{B}\mathbf{w}_t 
\end{align}
\end{subequations}
for $t \in \cT$, where $C_t: \mathbb{R}^S \mapsto \mathbb{R}$ is a time-varying smooth and convex function associated with the vector of states $\mathbf{y}_t$. Let $\bx_{t}^*$ be an optimal of~\eqref{eq:problem_DR}; the objective  is then to identify an optimal \emph{trajectory} $\{\bx_{t}^*, t \in \cT\}$. Before proceeding, a couple of examples of applications are provided. 
 
\vspace{.1cm} 
 
\emph{Example 1: Feeder-level problem}. For a feeder, $\mathbf{y}_t$ can collect voltages at some selected nodes~\cite{Dallanese2018feedback} and the net powers measured at the point of connection of the feeder with the rest of the grid. One may want to drive the state $\mathbf{y}_t$ towards a time-varying reference point $\mathbf{y}_{\text{ref},t}$ using the function $C_t(\mathbf{x}_t)=\frac{\beta}{2}\norm{ \mathbf{A} \mathbf{x}_t + \mathbf{B}\mathbf{w}_t - \mathbf{y}_{\text{ref},t}}^2$, with $\beta > 0$. In this case, $\mathbf{A}$ can be constructed based on the Jacobian of the power flow equations, linear approximations  of the  power flow equations, or by estimating the sensitivities of the network. As shown shortly, the proposed algorithmic framework does not need knowledge of the matrix $\mathbf{B}$.

\vspace{.1cm} 

\emph{Example 2: Neighborhood-level problem}. For an aggregations of devices in a neighborhood or community, $y_t$ represents the total active power at the point of interconnection of the rest of the grid. In this case, $\mathbf{A}$  boils down to a row-vector with all ones and $x_{m,t}$ represents the active power setpoints of the devices. In the spirit of a ``virtual power plant,'' $y_{\text{ref},t}$ can be a time-varying reference signal for the active power at the point of interconnection (to provide, for example, primal or secondary grid services). Section~\ref{sec:results} will illustrate this case. 
\vspace{.1cm}

However, solving problem \eqref{eq:problem_DR} at each time step $t$ might be not viable because of the main challenges \emph{C1)}-\emph{C2)}; more specifically, one may not be able to collect measurements of the non-controllable powers $\bw_t$ because of sensing limitations, and because the function $U_m$ may be unknown or largely different from synthetic models. In this paper, we propose a  feedback-based online algorithm where: \emph{i)} measurements of $\mathbf{y}_t$ are utilized to estimate the gradient of the function $C_t \left(\mathbf{A} \mathbf{x}_t + \mathbf{B}\mathbf{w}_t \right)$; and, \emph{ii)} feedback from the users are utilized to estimate the functions $\{U_m\}_{m = 1}^M$ concurrently with the execution of the online algorithm.     
In this paper, the function $U_m$ is estimated using feedback information from the user via GPs. Specifically, a shape-constrained GP approach~\cite{jrnl_wang} is pursued to approximate the discomfort function with a strongly convex and smooth function. 
Accordingly, let $\hat{U}_{m}(x_{m,t})$ be the estimate of $U_{m}(x_{m,t})$ available at time $t$. In lieu of~\eqref{eq:problem_DR}, the goal is then to identify solutions of the following optimization problem in an online fashion: 
\begin{equation}
    \mathbf{x}_{t}^* = \underset{\{x_{m,t} \in \mathcal{X}_{m,t}\}_{i=1}^M}{\text{argmin}} \; \sum_{m=1}^M \hat{U}_{m}(x_{m,t}) + C_t \left(\mathbf{A} \mathbf{x}_t + \mathbf{B}\mathbf{w}_t \right).
    \label{eq:problem_DR2}
\end{equation}

How to construct $\hat{U}_{m}(x_{m,t}) $ is explained  next.

\subsection{Shape-constrained Gaussian Processes}
\label{sec:shapeGP}

In this section, we introduce the main concepts underpinning GPs~\cite{GPML_Rasmussen06} and shape-constrained GPs~\cite{jrnl_wang}. They both offer a non-parametric model that is convenient for the learning setting of this paper because of the simplicity of the online updates and the ability to handle asynchronous and noisy data. In this section, the subscripts $m$ and $t$ are removed under the understanding that the technical arguments apply to each of the discomfort functions for all times. 

\subsubsection{Gaussian Process in a nutshell} \label{secGP}

A GP is a stochastic process $U(x)$ and it is specified by its mean function $\mu({x})$ and its covariance function $k(x, x')$;
i.e., for any $x, x' \in \mathcal{X} \subseteq \mathbb{R}$, 
$\mu({x}) = \mathbb{E}[U({x})]$ and $k({x}, {x}') = \mathbb{E}[(U({x}) - \mu({x}))(U({x}') - \mu({x}'))]$ \cite{GPML_Rasmussen06}. Let $\mathbf{x}_p = [ x_1 \in \mathcal{X}, \dots, x_p \in \mathcal{X}]^\top$ be the set of $p$ sample points and let ${z}_i = U({x}_i) + \epsilon_i$, with $\epsilon_i \stackrel{iid}{\sim} \mathcal{N}({0}, \sigma^2)$ Gaussian noise, be the noisy measurements at the sample points $x_i \, \forall \, i=1, \dots, p$; and, define $\mathbf{z}_p = [z_1, \dots, z_p]^\top$. Then, the posterior distribution of $(U(x)|\mathbf{x}_p, \mathbf{z}_p)$ is a GP  with mean $\mu_p({x})$, covariance $k_p({x}, {x}')$, and variance $\sigma_p^2(x)$ given by:
\begin{subequations}
\begin{align}
    \mu_p({x}) &= \mathbf{k}_p(x)^\top (\mathbf{K}_p + \sigma^2 \mathbf{I}_p)^{-1} \mathbf{z}_p \\
    k_p({x}, {x}') &= k(x,x') - \mathbf{k}_p(x)^\top (\mathbf{K}_p + \sigma^2 \mathbf{I}_p)^{-1} \mathbf{k}_p(x') \\
    \sigma_p^2(x) &= k_p(x,x)
\end{align}
\end{subequations}
where $\mathbf{k}_p(x) = [k(x_1,x), \dots, k(x_p,x)]^\top$, $\mathbf{K}_p$ is the positive definite kernel matrix $[k(x,x')]$, and the subscipt $p$ indicates the number of data points in $\mathbf{x}_p$. Thus, an estimate of the (unknown) function $U({x})$ can be written as $U({x}) \sim \mathcal{GP}( \mu_p({x}),  k_p({x}, {x}'))$. The covariance function specifies the covariance $\text{cov}(U({x}),U({x}')) $ between pairs of random variables; using squared exponential (SE) kernel as an example,  it is  defined as
\begin{equation}
    k({x},{x}') = \sigma_{f}^2 \; e^{-\frac{1}{2l^2}({x}-{x}')^{2}}
\end{equation}
\noindent
for the univariante input case, 
where the hyperparameters are the  variance $\sigma_{f}^2$ and the characteristic length-scale $l$.

\subsubsection{The Derivative Processes of GP} 

It is convenient to consider the SE covariance function because the resulting process has derivatives of all orders (see, e.g., \cite[Theorem 2.2.2]{Adler_RF}). Since differentiation is a linear operator, derivatives of the GP remains a GP \cite{GPML_Rasmussen06}. To obtain a strongly convex function, we will use the second derivative process of the GP. In particular, the corresponding mean and covariance function (jointly with the original process and the second-order derivative process) are \cite{jrnl_wang}:

\begin{subequations}
\label{eq:GP_dsecond}
\begin{equation}
   \hspace{-4.0cm} \mathbb{E}\left[\frac{\partial^2 U({x})}{\partial {x}^2}\right] = \frac{\partial^2 \mu({x})}{\partial {x}^2} = {0}
\end{equation}
\vspace{-0.5cm}
\begin{multline}
        k^{22}(x,x') := \text{cov}\left[\frac{\partial^2 U({x})}{\partial {x}^2}, \frac{\partial^2 U({x}')}{\partial {x}'^2}\right] = \sigma_{f}^2  e^{-\frac{1}{2l^2}(x-x')^{2} } \times \\ \frac{1}{l^4}\left(\frac{1}{l^4}(x-x')^4 - \frac{1}{l^2}6(x-x')^2+3\right) 
\end{multline}
\vspace{-0.5cm}
\begin{multline}
k^{02}(x,x') :=    \text{cov}\left[\frac{\partial^2 U({x})}{\partial {x}^2}, U({x}') \right]  \\ = \sigma_{f}^2 e^{ -\frac{1}{2l^2}(x-x')^{2} }  \left(\frac{1}{l^4}(x-x')^2 - \frac{1}{l^2}\right) 
\end{multline}
\end{subequations}

\subsubsection{Shape Constraints}

Suppose that one acquires noisy observations $\bz_p$ of the GP at $p$ points $\bx_p$ (based on, e.g., the user's feedback), but no observations over the derivative process are available. However,   we will impose derivative constraints at $q$ points $\mathbf{d} := [d_1, \dots, d_q]^\top$ \cite{jrnl_wang}; that is, constraints on the shape of the function are imposed even at points where there is not observation of the actual process.  

Let $\mathbf{U}(\mathbf{x})= [U(x_1), \dots, U(x_p)]^\top$ and $\mathbf{U}''(\mathbf{d})= [U''(d_1), \dots, U''(d_q)]^\top$; then, the joint distribution of the GP and its second-order derivative is:
\begin{equation*}
    \begin{bmatrix} 
    \mathbf{U}(\mathbf{x}) \\
    \mathbf{U}''(\mathbf{d})
    \end{bmatrix}
    \sim \mathcal{N} \bigg( \begin{bmatrix} 
    \mu \mathbf{1}_p\\
    \mathbf{0}_q
    \end{bmatrix},
    \begin{bmatrix} 
    \mathbf{K}(\mathbf{x},\mathbf{x}) & \mathbf{K}^{\mathbf{02}}(\mathbf{x},\mathbf{d}) \\
    \mathbf{K}^{\mathbf{20}}(\mathbf{d},\mathbf{x}) & \mathbf{K}^{\mathbf{22}}(\mathbf{d},\mathbf{d})
    \end{bmatrix} \bigg),
\end{equation*}
where $\mathbf{K}(\mathbf{x},\mathbf{x}) = \mathbf{K}_p$, $\mathbf{K}^{\mathbf{02}}(\mathbf{x},\mathbf{d}) = [k^{02}(x, d)]$, $\mathbf{K}^{\mathbf{20}}(\mathbf{d},\mathbf{x}) = \mathbf{K}^{\mathbf{02}}(\mathbf{x},\mathbf{d})^\top$ and $\mathbf{K}^{\mathbf{22}}(\mathbf{d},\mathbf{d}) = [k^{22}(d, d')]$.


In the following, we will impose  constraints  via  indicator  functions. Assign to $U(\cdot)$ a GP prior, and consider obtaining an estimated function that is $L_U$-smooth and $\gamma_U$-strongly convex, for a given $L_U > 0$ and $\gamma_U > 0$. We adapt the results presented in \cite{jrnl_wang} for the marginal constrained prior distribution.


Following \cite[Lemma 3.1]{jrnl_wang}, the joint conditional posterior distribution of $(\mathbf{U}({x}^\circ)|\mathbf{U}''(\mathbf{d}), \mathbf{x}_p, \mathbf{z}_p)$, for a point ${x}^\circ$ of a new set of $p^\circ$ points, given the current observations $\mathbf{z}_p$, is a GP with mean, covariance, and standard deviation given by: 
\begin{subequations}
\label{eq:shape_GP}
    \begin{align}
        \label{eq:GP_def_mu}
        & \bar{\mu}_{p^\circ}({x^\circ}) = \mu \mathbf{1}_{p^\circ} + B_3(\bx, {x}^\circ,\mathbf{d}) B_1(\bx,\mathbf{d})^{-1}(\mathbf{z}_p - \mu \mathbf{1}_p) \nonumber \\
        & + (A_2({x}^\circ,\mathbf{d}) - B_3(\bx, {x}^\circ,\mathbf{d}) B_1(\bx,\mathbf{d})^{-1} A_1(\bx,\mathbf{d}))\mathbf{U}''(\mathbf{d}), \\
        & \bar{k}_{p^\circ}({x^\circ}, {x^\circ}') = A(\bx, {x}^\circ,\mathbf{d}), \\
        \label{eq:GP_def_sigma}
        & \bar{\sigma}_{p^\circ}(x^\circ) = \sqrt{A(\bx, {x}^\circ,\mathbf{d})},
    \end{align}
\end{subequations}
and the posterior distribution of $(\mathbf{U}''(\mathbf{d})| \mathbf{x}_p, \mathbf{z}_p)$ is given by:
\begin{equation*}
    (\mathbf{U}''(\mathbf{d})|\mathbf{x}_p,\mathbf{z}_p) \propto \mathcal{N}(\bmu(\mathbf{d}), \mathbf{D}(\mathbf{d}, \mathbf{d})) \mathbf{1}_{\{ \gamma_U \leq U''(d_i) \leq  L_U, \, i=1, \dots, q \}}
\end{equation*}
where $(\mathbf{U}''(\mathbf{d})|\mathbf{x}_p,\mathbf{z}_p)$ is a truncated normal distribution and, \vspace{-0.5cm}
\begin{subequations}
    \begin{align*}
        \bmu(\mathbf{d}) &= \mathbf{K}^{\mathbf{20}}(\mathbf{d},\mathbf{x})(\sigma^2 \mathbf{I} + \mathbf{K}(\mathbf{x},\mathbf{x}))^{-1}(\mathbf{z}_p - \mu \mathbf{1}_p),\\
        \mathbf{D}(\mathbf{d}, \mathbf{d}) &= \mathbf{K}^{\mathbf{22}}(\mathbf{d},\mathbf{d}) \\
        & \hspace{.5cm} -\mathbf{K}^{\mathbf{20}}(\mathbf{d},\mathbf{x})(\sigma^2 \mathbf{I} + \mathbf{K}(\mathbf{x},\mathbf{x}))^{-1}\mathbf{K}^{\mathbf{02}}(\mathbf{x},\mathbf{d}),\\
        A_1(\bx, \bd) &= \mathbf{K}^{\mathbf{02}}(\mathbf{x},\mathbf{d})\mathbf{K}^{\mathbf{22}}(\mathbf{d},\mathbf{d})^{-1}, \\
        A_2({x}^\circ,\mathbf{d}) &= \mathbf{K}^{\mathbf{02}}({x}^\circ,\mathbf{d})\mathbf{K}^{\mathbf{22}}(\mathbf{d},\mathbf{d})^{-1},\\
        B_1(\mathbf{x},\mathbf{d}) &= \sigma^2 \mathbf{I} + \mathbf{K}(\mathbf{x},\mathbf{x}) \\
        & \hspace{.5cm} - \mathbf{K}^{\mathbf{02}}(\mathbf{x},\mathbf{d})\mathbf{K}^{\mathbf{22}}(\mathbf{d},\mathbf{d})^{-1}\mathbf{K}^{\mathbf{20}}(\mathbf{d},\mathbf{x})\\
        B_2({x}^\circ,\mathbf{d}) &= \mathbf{K}({x}^\circ,{x}^\circ) \\
        &\hspace{.5cm}  - \mathbf{K}^{\mathbf{02}}({x}^\circ,\mathbf{d})\mathbf{K}^{\mathbf{22}}(\mathbf{d},\mathbf{d})^{-1}\mathbf{K}^{\mathbf{20}}(\mathbf{d},{x}^\circ),\\
        B_3(\bx, {x}^\circ,\mathbf{d}) &= \mathbf{K}({x}^\circ,\mathbf{x}) \\
        & \hspace{.5cm} - \mathbf{K}^{\mathbf{02}}({x}^\circ,\mathbf{d})\mathbf{K}^{\mathbf{22}}(\mathbf{d},\mathbf{d})^{-1}\mathbf{K}^{\mathbf{20}}(\mathbf{d},\mathbf{x}) ,\\
        A(\bx,{x}^\circ,\mathbf{d})) &= B_2({x}^\circ,\mathbf{d}) \\
        & \hspace{.5cm} -  B_3(\bx,{x}^\circ,\mathbf{d}) B_1(\mathbf{x},\mathbf{d})^{-1} B_3(\bx,{x}^\circ,\mathbf{d})^\top
    \end{align*}
\end{subequations}
with $\mu$ and $\sigma^2$ given parameters of the prior. The parameters $l$ and $\sigma_f^2$ of the GP can be estimated, for example, by using the maximum likelihood estimator \cite{GPML_Rasmussen06}. The locations of the virtual derivative points are defined beforehand.  By imposing the smooth and strong convexity constraints on points $\mathbf{d}$ that are dense enough, shape-constrained GPs ensure that the posterior mean function $\bar{\mu}_{p^\circ}(x^{\circ})$ is ``practically'' (i.e., indistinguishable for all practical purposes) smooth and strongly convex~\cite{jrnl_wang}. The choice of shape-constrained GPs versus exact methods, such as smooth strong convex regression~\cite{SimonettoCR2020} (which would ensure shape properties exactly and everywhere) is motivated by the fact that the latter is more computationally intensive and its learning rate can be significantly slower.

\section{Online GP-based Demand Response}

When the functions $\{U_m\}_{m = 1}^M$ are known and the non-controllable powers can be  measured at each time instant $t$, then the time-varying problem~\eqref{eq:problem_DR2} can be solved in an online fashion using the following online projected gradient method:  
\begin{align}
\label{eq:online-gradient}
    \bx_{t} = \textrm{proj}_{\mathcal{X}_{t}} \{\bx_{t-1} - \alpha \left(\nabla_\bx U(\bx_{t-1}) + \nabla_\bx C_{t}(\bx_{t-1}) \right)\}
\end{align}
where $U(\bx) := \sum_{m = 1}^M U_m(x_m)$ for brevity, 
$\textrm{proj}_{\mathcal{X}}\{\by\} := \arg \min_{\bx \in \mathcal{X}} \|\bx - \by\|^2$ is the projection operator, and $\alpha > 0$ is the step size. To address the challenges \emph{C1)}-\emph{C2)}, the online algorithm~\eqref{eq:online-gradient} is modified as explained next.

\subsection{Online Algorithm}

Recall that $t \in \cT$ is the time index. We now introduce an additional index $p_m(t)$ (one per device or user), used as a counter for the number of data points $\bz_{m,t}:= [z_{m,1}, \ldots, z_{m,p_m(t)}]^\top$ received from the $m$th  user up to time $t$; we recall that $z_{m,p_m(t)} = U_m(x_{m,t}) + \epsilon_{m,t}$ ($t$ being the time when the $p_m$th user feedback is received). The counter $p_m(t)$ does not generally coincide with $t$, since a user may provide feedback sporadically or at a slower time scale (whereas the algorithm is run an a fast time scale). Hereafter, we omit the dependence of $p_m$ on $t$ for notation simplicity. 

With $p_m$ data points 
available (i.e., received from the $m$th user), we define the estimate $\hat{U}_{m,p_m}(x_{m,t})$ of $U_{m}(x_{m,t})$ as:
\begin{equation}
    \label{eq:estimate_U}
    \hat{U}_{m,p_m}(x_{m,t}) := \bar{\mu}_{m,p_m}(x_{m,t})
\end{equation}
where $\bar{\mu}_{m,p_m}(x_{m,t})$ 
is given by \eqref{eq:GP_def_mu} based on $p_m$ data points (which we remind that is ``practically'' smooth and strongly convex). In other words, $\hat{U}_{m,p_m}(x_{m,t})$ is obtained via the mean of the shape-constrained GP when feedback from the user is received $p_m$ times. Further, at a given point $x_{m,t}$, the derivative of $\hat{U}_{m,p_m}(x_{m,t})$ is estimated via finite-difference as~\cite{flaxman2004online}:
\begin{equation}
\label{eq:gradient_U}
    v_{m,p_m}(x_{m,t}) := \frac{\hat{U}_{m,p_m}(x_{m,t} + \delta) - \hat{U}_{m,p_m}(x_{m,t})}{\delta} 
\end{equation}
with $\delta$ a pre-selected parameter. For future developments, let $\mathbf{v}(\mathbf{x}_t) := [v_{1,p_1}(x_{1,t}), ..., v_{M,p_M}(x_{M,t})]^\top$.  

The evaluation of the gradient of $C(\mathbf{x}_t)$ requires measurements of the non-controllable devices $\mathbf{w}_t$ at each time step $t$. Similar to, e.g.,~\cite{bolognani2013distributed,Dallanese2018feedback},  measurements $\hat{\by}_t$ can be utilized in the computation of the gradient of $C(\mathbf{x}_t)$ instead of the map   $\mathbf{y}_t = \bA \bx_t + \bB \bw_t$. For example, if the function $C_t(\mathbf{x}_t)$ is $C_t(\mathbf{x}_t)=\frac{\beta}{2}\norm{ \mathbf{A} \mathbf{x}_t + \mathbf{B}\mathbf{w}_t - \mathbf{y}_{\text{ref},t}}^2$, its gradient reads $\nabla C_t(\bx_t) = \beta \mathbf{A}^\top (\bA \bx_t + \bB \bw_t - \mathbf{y}_{\text{ref},t})$; on the other hand, an estimate of the gradient using the measurement $\hat{\by}_t$ amounts to $\mathbf{s}_t := \beta \mathbf{A}^\top (\hat{\mathbf{y}}_t - \mathbf{y}_{\text{ref},t})$. Indeed, $\mathbf{s}_t$ can be interpreted as a noisy version of $\nabla_\mathbf{x} C(\mathbf{x}_t)$ \cite{Mag_signal_2020}.

Overall, the proposed shape-constrained GP-based online projected gradient descent (SGP-OPGD) method involves the sequential execution of the following step:
\begin{align}
\label{eq:online-gradient-GP}
    \bx_{t} = \textrm{proj}_{\mathcal{X}_{t}} \{\bx_{t-1} - \alpha \left(\bv(\bx_{t-1}) + \bs_{t} \right)\}
\end{align}
where we recall that $\bv(\bx_{t})$ is an estimate of the gradient of $\hat{U}(\bx_t)$, where $\hat{U}(\bx_t) := \sum_{m = 1}^M \hat{U}_{m,p_m}(x_{m,t})$, $\mathbf{s}_t$ is a noisy version of $\nabla_\mathbf{x} C(\mathbf{x}_t)$, $t$ represents the time index, and $p_m$ is the data counter for the user's feedback per device. 

The steps of the SGP-OPGD are detailed in Algorithm \ref{algo}. 
Notice that the update of $\bx_t$ decouples into $M$ parallel steps (one per device); this enables a distributed setting with a so-called ``gather-and-broadcast'' architecture where measurements of $\hat{\by}_t$ are collected at a central location, $\mathbf{s}_t$ is broadcasted to the devices, and $x_{m,t}$ is computed locally at each device. Further, the function $\hat{U}_{m,p_m}(x_{m,t})$ is computed locally. 

\begin{algorithm}[ht]

\textbf{Initialize}: $\bx_0$,  $\alpha = \frac{2}{\gamma + L}$; prior on $\{\hat{U}_m\}_{m = 1}^M$ if available. \\

\nl \bf for $t = 1, 2, \dots, T $ \bf do \\
\nl \quad \quad \normalfont{Collect measurement} $\hat{\mathbf{y}}_t$ \\ 

\nl \quad \quad Compute the estimate gradient  $\mathbf{s}_{t}$ \\

\nl \quad \bf for $m = 1, 2, \dots M$ \bf do \\
\nl \quad \quad \bf if \normalfont{Feedback is given: } \\ 
\nl \quad \quad \quad $ p_m \rightarrow p_m+1$\\
\nl \quad \quad \quad \normalfont{Collect} $z_{m,p_m}$ \normalfont{and add it to} $\bz_{m,t}$ \\ 
\nl \quad \quad \quad \normalfont{Update} $\hat{U}_{m,p_m}(x_{m,t})$ \normalfont{and compute} $v_{m,p_m}(x_{m,t})$ \\
\nl \quad \quad \bf else \normalfont{Keep} $\hat{U}_{m,p_m}(x_{m,t-1})$ and $v_{m,p_m}(x_{m,t-1})$\\
\nl \quad \quad \normalfont{Update setpoint as} \\ 
\vspace{-0.5cm}
\begin{equation*}
    \qquad x_{m,t} = \textrm{proj}_{\mathcal{X}_{m}} \{x_{m,t-1} - \alpha  (v_{m,p_m}(x_{m,t-1}) + s_{m,t})\}
\end{equation*}
\nl \quad \bf end for \\
\nl \bf end for
    \caption{{SGP-OPGD method} 
    \label{algo}}
    
\end{algorithm}

\subsection{Analysis}

The convergence of the online algorithm is compared against the optimal trajectory $\{\bx^*_t\}_{t \in \cT}$ and the optimal value function of~\eqref{eq:problem_DR2}. Presuming a synthetic function $U(\bx)$, the difference between $\bx^*_t$ and a solution of~\eqref{eq:problem_DR} will be assessed numerically in Section~\ref{sec:results}. Hereafter, we define $f_t(\bx) :=  \hat{U}(\bx) + C_t(\bx)$ 
for brevity.   

We begin with the following standard assumptions.
\begin{enumerate}
    \item [\textit{AS1:}] The function $f_t$ is $L$-smooth on $\mathcal{X}$; i.e., $\norm{\nabla f_t(\mathbf{x})-\nabla f_t(\mathbf{x}')} \leq L \norm{\mathbf{x} - \mathbf{x}'}$ for all $t \in \mathcal{T}$ and $\mathbf{x}$, $\mathbf{x'} \in \mathcal{X}$.

    
    \item[\textit{AS2:}] The function $f_t$ is $\gamma$-strongly convex. 
    
    \item[\textit{AS3:}] The inexact gradient $\tilde \nabla f_t(\mathbf{x}) = \mathbf{v}_p(\mathbf{x}) + \mathbf{s}_t$ is defined as $\tilde \nabla f_t(\mathbf{x}_t) :=  \nabla f_t(\mathbf{x}_t) + \mathbf{e}_{1,t} + \mathbf{e}_{2,t}$, where $\mathbf{e}_{1,t}$ is the error in the gradient of $\{ \hat{U}_{m,p_m}\}_{m=1}^M$ and $\mathbf{e}_{2,t}$ is the error in the estimated gradient $\mathbf{s}_t$. The sequence $\{ \mathbf{e}_t := \mathbf{e}_{1,t} + \mathbf{e}_{2,t} \in \mathbb{R}^M, t \in \mathcal{T} \}$ is bounded; i.e., $\norm{\mathbf{e}_t} < \infty$. 
    

\end{enumerate}

Regarding \textit{AS1}, $L$ is given by $L = L_U + L_C$, with $L_U$ and $L_C$ the Lipshitz constants of the gradients of $\hat{U}$ and $C_t$, respectively; notice that the Lipshitz constant of the gradient of each individual function $\hat{U}_{m,p_m}$ is set a priori as in~\eqref{eq:shape_GP}. If $C_t$ is convex but not strongly convex, only the strong convexity coefficient of $\hat{U}$ plays a role in \textit{AS2} [cf.~\eqref{eq:shape_GP}].   

The variation between any two consecutive optimal points is defined as $r_t := \norm{\mathbf{x}_{t-1}^* - \mathbf{x}_{t}^*}$.
Now, define the path length and the cumulative  gradient error as~\cite{Mag_signal_2020,Jadbabaie2015}
\begin{equation}
    \omega_T := \sum_{t=1}^{T} r_t, \qquad  E_T = \sum_{t=1}^{T} \norm{\mathbf{e}_t}.
    \label{eq:para1}
\end{equation}
These metrics will be utilized in the following results. 

\vspace{.1cm} 

\begin{proposition}
Assume that $\alpha \in (0, 2/L)$. Under Assumptions \textit{AS1}-\textit{AS3}, the SGP-OPGD algorithm constructs a sequence $\{ \mathbf{x}_{t} \}_{t \in \mathcal{T}}$ such that
\begin{equation}
    \norm{\mathbf{x}_{t} - \mathbf{x}_t^*}
    \leq \rho \norm{\mathbf{x}_{t-1} - \mathbf{x}_{t-1}^*} + \rho r_t + \alpha \norm{\mathbf{e}_t},
    \label{eq:th1}
\end{equation}
where $\rho := \text{max}\{|1-\alpha \gamma|, |1 - \alpha L| \} < 1$.
\label{th1}
\end{proposition}

\vspace{.1cm} 

\begin{corollary}
Under assumption \textit{AS1}-\textit{AS3}, with $\alpha \in (0, 2/L)$, the cumulative tracking error of the SGP-OPGD algorithm can be bounded as:
\begin{equation}
    \sum_{t=1}^T \norm{\mathbf{x}_{t} - \mathbf{x}_t^*} \leq \frac{1}{1-\rho} \left[ \rho \norm{\mathbf{x}_0 - \mathbf{x}_{0}^*} + \rho \omega_T + \alpha E_T \right].
    \label{eq:th2}
\end{equation}
\label{th2}
\end{corollary}

\vspace{.1cm} 

Proposition \ref{th1} establishes $Q$-linear convergence to a bounded error of the SGP-OPDG algorithm~\cite{Mag_signal_2020}; that is, each step of the algorithm is contractive up an error $\alpha \norm{\mathbf{e}_t} + r_t$ given by the temporal variability of the problem and the errors in the gradient computation. On the other hand, Corollary \ref{th2} asserts that the tracking error of the algorithm is bounded if $E_T$ and $\omega_T$ grow as $\mathcal{O}(T)$, and it goes to zero asymptotically if $E_T$ and $\omega_T$ grow sublinearly  in $T$; that is, if they grow as $o(T)$.  

Finally, we provide a bound on the dynamic regret next.

\vspace{.1cm}


\begin{proposition} 
\label{thm:convergence}
Suppose that Assumptions \textit{AS1}-\textit{AS3} hold, and let $\alpha \in (0, 2/L)$. Then, the dynamic regret of  the SGP-OPGD algorithm can be bounded as:
\begin{equation}
\frac{1}{T} \sum_{t=1}^T [f_t(\mathbf{x}_t) - f_t(\mathbf{x}_t^*)] = \mathcal{O}\left(T^{-1}\omega_T + T^{-1} E_T \right).
\end{equation}
\end{proposition} 

\vspace{.1cm}

In par with Corollary \ref{th2}, the dynamic regret  is sublinear if $\omega_T$ and $E_T$ are both sublinear. If $\omega_T$ and $E_T$ grow linearly, then the dynamic regret behaves as $\mathcal{O}$(1). Moreover, if the gradient is obtained without error (i.e., $\mathbf{e}_t = \mathbf{0}$), the dynamic regret exhibits a behavior similar to \cite{IP-OGD}.

The proofs follow steps similar to \cite{Inexact_prox_online,IP-OGD}; Proposition~\ref{thm:convergence} uses the fact that, from the  continuity  of  the  gradient  and  the compactness  of $\cX_t$, the norm of the gradient is bounded. The proofs are not provided here due to the page limit. 


\section{Illustrative Results}
\label{sec:results}

We consider a neighborhood-level problem as in  \textit{Example 2} in Section \ref{sec:modeling}. In the considered example, we control 15 batteries, 10  HVAC units (equipped with variable speed drives), and 5 electric vehicles (EVs). The objective is to maintain the aggregate active power $\sum_m x_{m,t}$ close to a reference point $y_{\text{ref},t}$ while minimizing the discomfort/dissatisfaction for each user. The operational sets for the devices are: \textit{(i)} batteries constraints $\mathcal{X}_m = [-8,8] \text{ kW } \forall \, m = 1, \dots, 15$; \textit{(ii)} HVAC constraints $\mathcal{X}_m = [5, 15] \text{ kW } \forall \, m = 16, \dots, 25$; and \textit{(iii)} EV constraints $\mathcal{X}_m = [7, 50] \text{ kW } \forall \, m = 26, \dots, 30$. To concretely assess the performance of the shape-constrained GP, the discomfort functions $\{U_m\}_{m=1}^M$ are assumed to be quadratic; the minimum of each of the functions is inside the set constraints $\mathcal{X}_m$, and it corresponds to a preferred setting of the user. For example, for EVs they represent a preferred charging rate; for  HVAC systems, they represent a preferred temperature setpoint (converted into a preferred power setpoint)~\cite{Lesage_2020}. The function $C_t(\bx)$ is $C_t(\mathbf{x}_t)=\frac{\beta}{2}(\sum_m x_{m,t} + \mathbf{1}^\top \mathbf{w}_t - y_{\text{ref},t})^2$,  where the non-controllable loads are taken from the Anatolia dataset (National Renewable Energy
Laboratory, Tech. Rep. NREL/TP-5500-56610) and have a granularity of 1 second. 

\begin{figure}[!ht]
  \centering
  \subfigure[]{\includegraphics[width=0.5\textwidth]{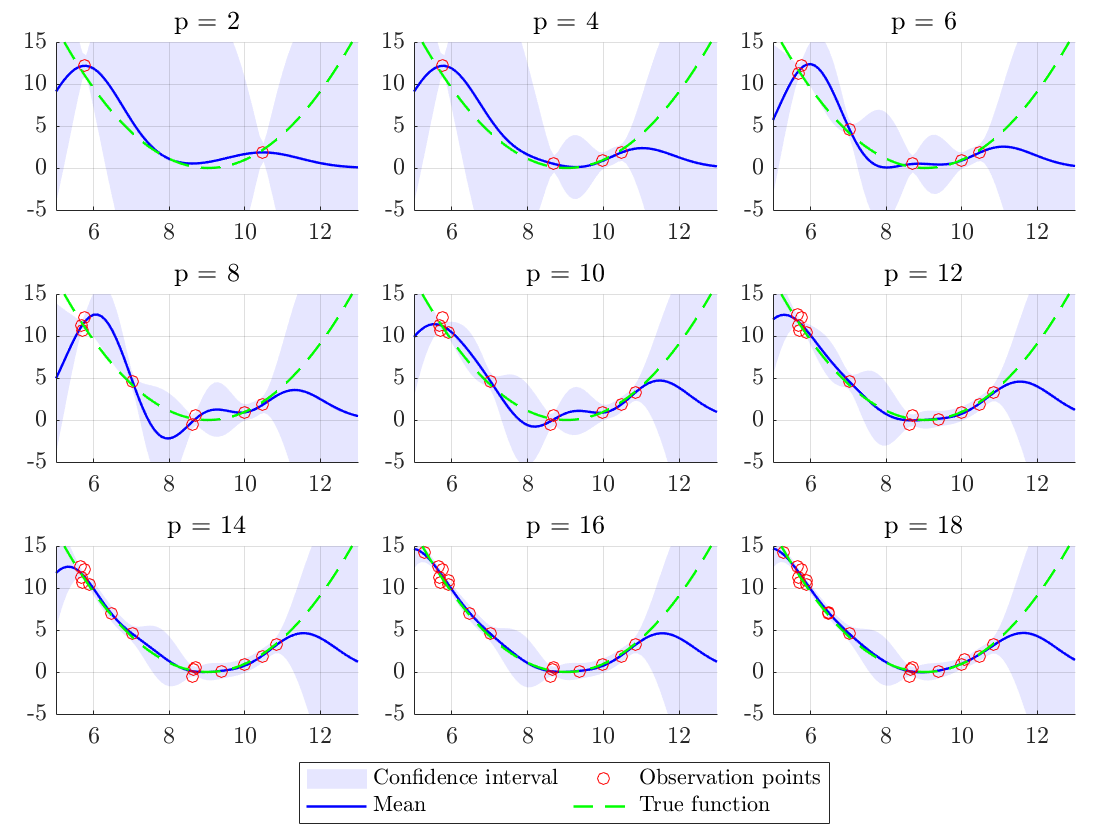}}
  \subfigure[]{\includegraphics[width=0.5\textwidth]{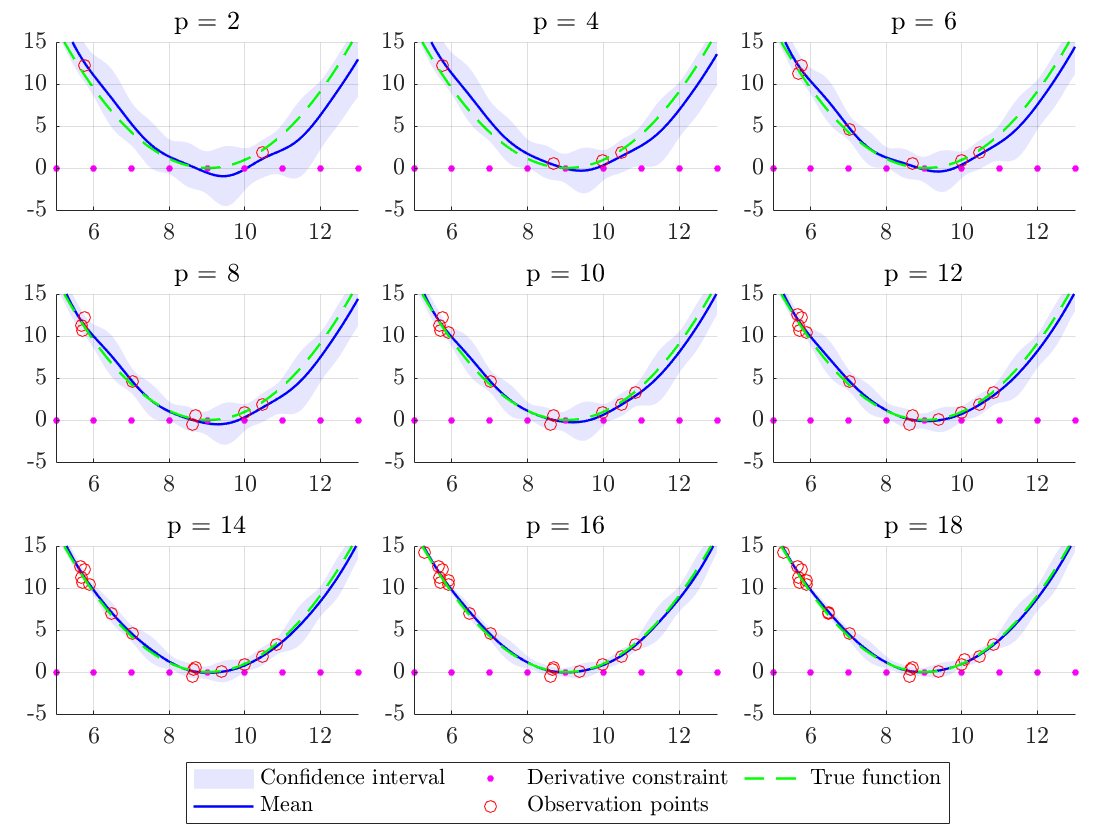}}
  \caption{Example of the estimation of the discomfort function $\hat{U}_{m,p_m}$  with hyperparameters $\sigma_f = 1$ and $l = 10$. (a) Standard GP regression and (b) shape-constrained GP regression. The observations points are a mixed of prior points and a sub-sequence produced by the online algorithm. One can see that the function estimated with shape-constrained GPs is ``practically'' smooth and strongly convex, as desired, after only a few feedback points $p_m$.}
  \label{fig:regular_shape_GP}
  \vspace{-.3cm}
\end{figure}

As an example of estimation of the discomfort functions using the shape-constrained GP, Figure \ref{fig:regular_shape_GP} illustrates the estimated function for a device for a different number of observations $p$; in particular, the estimated functions using a standard GP regression and the shape-constrained GP are illustrated. 

We run the online algorithm for a period of 12 hours staring at 12:00 am; each step of Algorithm \ref{algo} is  performed every 5 seconds (expect for HVAC, which are updated at a slower rate). A prior $\{ \hat{U}_{m,p_m} \}_{m=1}^M$ is determined from some  noisy measurements ($\sigma = 5$) and $\hat{U}_{m,p_m}(x_{m,t})$ is updated through user's feedback every 30 min. 
\begin{figure*}[!ht]
  \centering
  \includegraphics[width=1\textwidth]{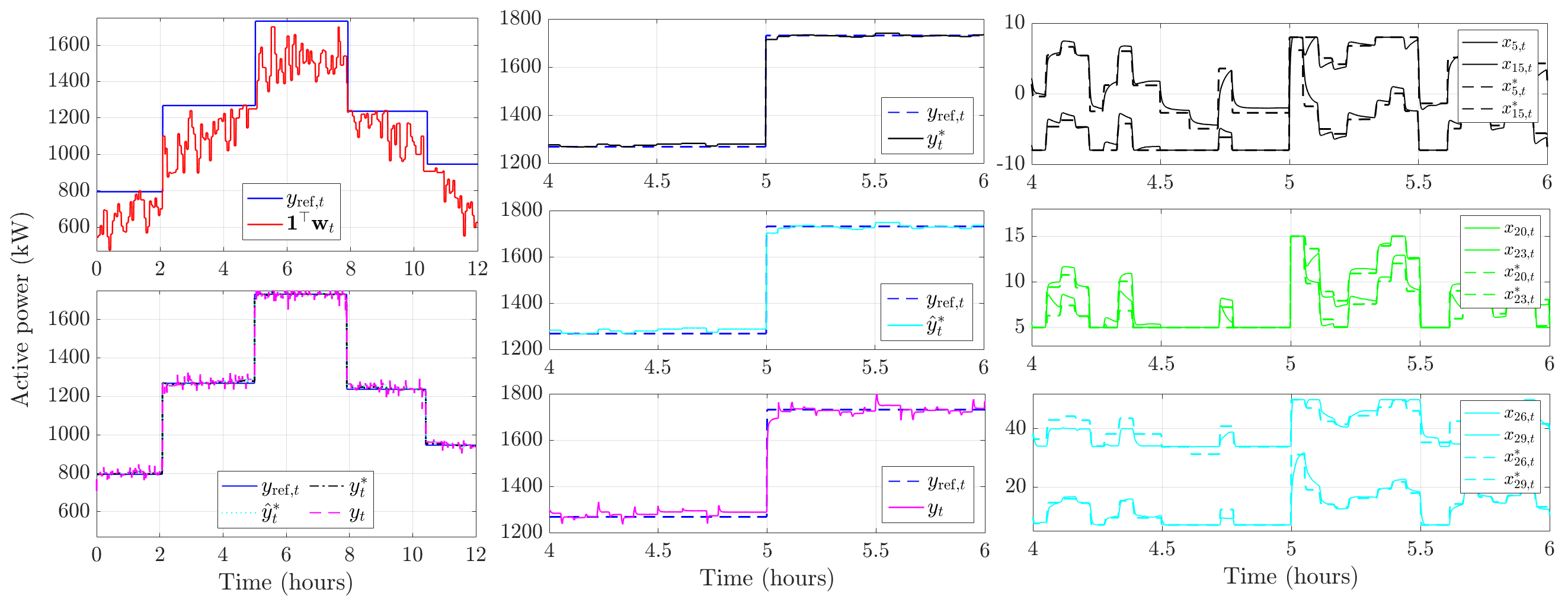}
  \caption{Solution of the SGP-OPGD algorithm. \emph{Left top}: references setpoints $y_{\text{ref},t}$ and aggregate non-controllable loads $\mathbf{1}^\top \bw_t$. \emph{Left bottom}: optimal trajectories $y^*_t$ for a known synthetic functions $\{ U_m \}_{m=1}^M$, learned trajectory $\hat{y}^*_t$ where $\hat{U}_{m,p_m}$ is used and  the problem is solved to convergence at each time, and  trajectory $y_t$ for the SGP-OPGD algorithm. \emph{Center}: Zoomed view for the trajectories ${y}^*_t$, $\hat{y}^*_t$ and ${y}_t$ on the time period 4:00 pm - 6:00 pm. \emph{Right}: example of active power setpoints for the SGP-OPGD method of 6 representative devices  on the time period 4:00 pm - 6:00 pm.
  %
  %
  }
  \label{fig:sim_case}
\end{figure*}
The results in  Figure \ref{fig:sim_case} for the the SGP-OPGD algorithm are compared with two trajectories: \textit{(i)} trajectory for the optimal solution $\mathbf{x}^*_t$ for a known synthetic discomfort functions $\{ U_m \}_{m=1}^M$, where the problem is solved to convergence; 
\textit{(ii)} trajectory for the learned optimal solution $\mathbf{\hat{x}}^*_t$ when  $\{ \hat{U}_{m,p_m} \}_{m=1}^M$ is estimated as in \eqref{eq:estimate_U}, where also the problem is solved to convergence. In this case, the estimate of the gradients for $\{ \hat{U}_{m,p_m} \}_{m=1}^M$ are calculated using a finite difference method; 21 noisy observations ($\sigma = 0.5$) for each $\{ U_m \}_{m=1}^M$ are used. For the online algorithm, the step-size is $\alpha = 0.002$, that corresponds to the optimal step-size for the online gradient descent algorithm. 

\begin{figure}[!t]
  \centering
  \includegraphics[width=0.5\textwidth]{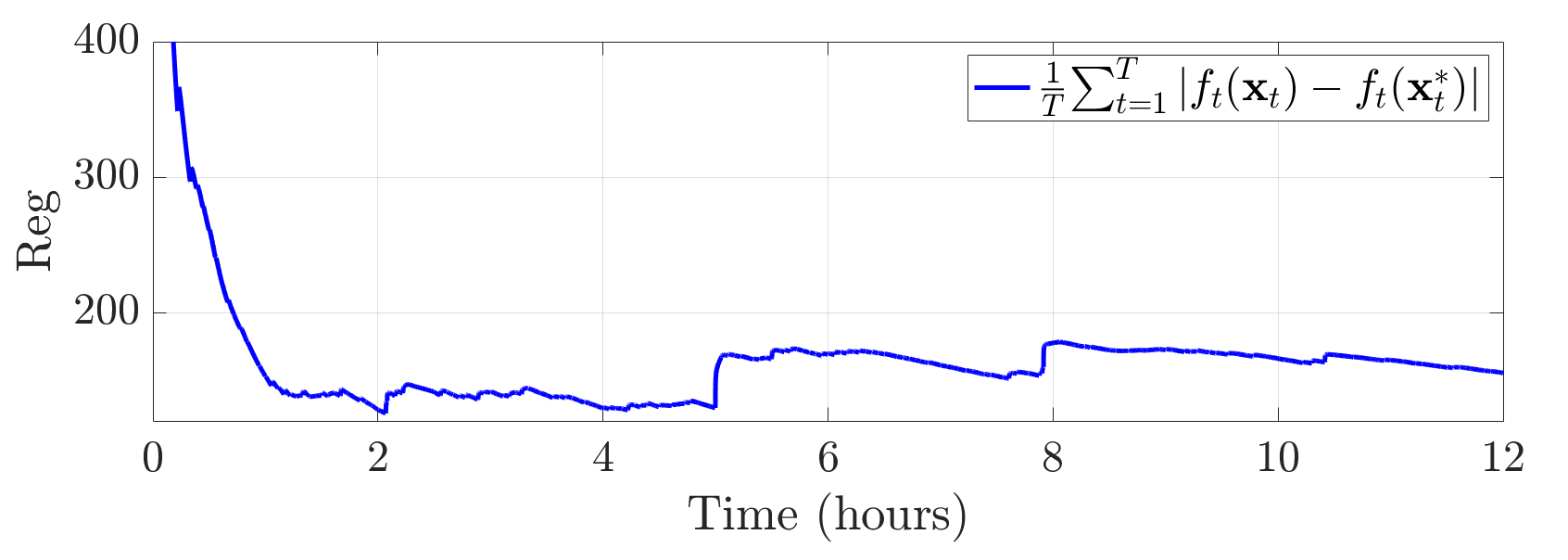}
  \caption{Dynamic regret of the SGP-OPGD algorithm. 
  }
  \label{fig:reg}
  \vspace{-.2cm}
\end{figure}

Figure \ref{fig:reg} shows the behavior of the performance metric for the SGP-OPGD algorithm, i.e., the  dynamic regret $\frac{1}{T} \sum_{T = 1}^T |f_t(\mathbf{x}_{t}) - f_t(\mathbf{x}_{t}^*)|$. It can be seen that the dynamic regret exhibits a $\mathcal{O}(1)$ asymptotic behavior; the jumps in the dynamic regret corresponds to instants where the reference $y_{\text{ref},t}$ changes abruptly.  

\bibliographystyle{IEEEtran}
\bibliography{References}

\end{document}